\documentclass [12pt, draft] {amsart}
\usepackage {amssymb, amsmath, amscd}

\theoremstyle {plain}
\newtheorem {theorem} {Theorem}
\newtheorem {corollary} {Corollary}

\newtheorem {lemma} {Lemma}

\theoremstyle {definition}
\newtheorem {definition} {Definition}

\theoremstyle {remark}

\newtheorem* {acknowledgments} {Acknowledgments}

\numberwithin {equation} {section}

\begin {document}
\title {Some endomorphisms of $II_1$ factors: Part II}
\author {Hsiang-Ping Huang}
\address {Department of Mathematics, University of Utah, Salt Lake City, UT 84112}
\email {hphuang@math.utah.edu}
\thanks {Research partially supported by National Center for Theoretical
         Sciences, Mathematics Division, Taiwan}

\pagestyle{plain}

\begin {abstract}
For any finite dimensional $C^*$-algebra $A$ with a trace vector
$\vec s$ whose entries are rational numbers, we give an endomorphism
$\Phi$ of the hyperfinite $II_1$ factor $R$  such that:
\[
  \forall \ k \ \in \ \mathbb {N}, \
  \Phi^k (R)' \cap R= \otimes^k A.
\]
The canonical trace $\tau$ on $R$ extends the trace vector $\vec s$
on $A$. Therefore the minimal projection is not necessarily equivalent
to each other.

\end {abstract}
\maketitle
\tableofcontents
\section {Introduction} \label {S: intro}

The study of subfactor theory is centered on describing
the position of a subfactor $N$ embedded into an ambient factor $M$.
The standard invariant 
associated with Jones basic construction, 
\[
  N \subset M \subset M_1 \subset M_2 \subset \cdots
\]
is a complete invariant in the amenable case.
To classify the standard invariant is 
the most important task ever since the birth of subfactor theory.
Many ground-breaking works have been done. Yet even more
puzzles remain unsolved.

For a hyperfinite $II_1$ subfactor of finite Jones index, 
it is equipped with an extra structure: an endomorphism $\Phi$,
sending the ambient factor $M$ onto the subfactor $N$.
Therefore it is only natural to investigate the "action",
mimicking A.Connes' marvellous work.

A well-known example is the canonical shift in a strongly 
amenable inclusion.
Another surprising example is  the binary shift \cite {rP88}
which gives rise to a counterexample that fails the tensor product formula
for entropy \cite{hN95}.
Via the Cuntz algebra, a lot of endomorphisms have been manufactured.

A nice result by M.Choda \cite {mC871}
states that $\Phi$ is outer-conjugate to
$\Psi$ as an endomorphism of $M$ onto $N$ if and only if $N_{\Phi}$
is conjugate to $N_{\Psi}$ as a subfactor of $M \otimes M_2({\mathbb C})$, 
where
\[
  N_\Phi= \{ \begin{bmatrix}x&0\\0&\Phi(x) \end{bmatrix} \ | \ 
  x \in M \} \ {\rm and \ } N_\Psi= \{ \begin{bmatrix}x&0\\0&\Psi(x) \end{bmatrix} 
  \ | \ x \in M \}
\]
Obviously we are interested in the Jones basic construction
\[
  N_\Phi \subset M \otimes M_2({\mathbb C})= \tilde M \subset {\tilde M}_1
  \subset {\tilde M}_2 \subset \cdots
\]

Unfortunately the standard invariant is hard to compute except for a 
basic endomorphism, where $\Phi$ can be extended to $M_n$ the tower algebra. 
Yet we do know the standard invariant contains a tower of finite dimensional 
$C^*$-algebras, the relative commutant algebras,
\[
  \{ \Phi^k (M)' \cap M \}_{k \in {\mathbb Z}_+} \subset 
  \{ \Phi^k (M)' \cap M_n \}_{k, n \in {\mathbb Z}_+}
\] 
For the moment, we are concentrated on the first tower,
which is easier to tackle with.

In \cite {hH} we coupled the notion of n-unitary shift with the
shift on $\otimes^{\infty}_{i=1} A$, for any finite dimensional
$C^*$-algebra $A$. We constructed an endomorphism $\Phi$ on $R$, 
which gives
\[
  \Phi^k (R)' \cap R= \otimes^k A.
\]
The Jones index $[R: \Phi (R)] = ({\rm  rank\ } (A))^2$, here ${\rm
rank\ } (A)$ is the dimension of the maximal abelian subalgebra of
$A$. The minimal projection of $A$ is equivalent to each other.

The major part of this paper is to generalize the above result 
to an arbitrary trace vector of rational entries on $A$.
Thus the minimal projection of $A$ may or may not be equivalent to
each other.

\section {Preliminaries} \label {S: prelim}

Let $M$ be a $II_1$ factor with the canonical trace $\tau$. Denote
the set of unital *-endomorphisms of $M$ by $End (M, \tau)$. Then
$\Phi \in End (M, \tau)$ preserves the trace and $\Phi$ is
injective. $\Phi (M)$ is a subfactor of $M$. If there exists a
$\sigma \in Aut (M)$ with $\Phi_1 \cdot \sigma= \sigma \cdot \Phi_2$
for $\Phi_i \in End (M, \tau)$ $(i= 1, 2)$ then $\Phi_1$ and
$\Phi_2$ are said to be conjugate. If there exists a $\sigma \in Aut
(M)$ and a unitary $u \in M$ such that $Ad u \cdot \Phi_1 \cdot
\sigma= \sigma \cdot \Phi_2$, then $\Phi_1$ and $\Phi_2$ are outer
conjugate.

The Jones index $[M: \Phi (M)]$ is an outer-conjugacy invariant. We
consider only the finite index case unless otherwise stated. In such
case, there is a distinguished outer-conjugacy invariant: the tower
of inclusions of finite dimensional $C^*$ algebras, $\{ A_k= \Phi^k
(M)' \cap M \}_{k=1}^{\infty}$.

\begin {lemma}
$A_k = \Phi^k (M)' \cap M$ contains an subalgebra that is isomorphic
to $\otimes_{l= 1}^k A_1$, the $k$-th tensor power of $A_1$, where
$A_1 = \Phi (M)' \cap M$ as denoted.
\end {lemma}

The dimension of the relative commutant ${\Phi^k (M)}' \cap M$ is
known to be bounded above by the Jones index $[M: \Phi(M)]^k$. Lemma
1 provides the lower bound for the growth estimate.

\begin {lemma}
For any finite dimensional $C^*$-algebra $A$, there exists a $\Phi
\in End (R, \tau)$ such that the relative commutant $ \Phi^k (R)'
\cap R$ is isomorphic to $\otimes_{i=1}^k A$. Here $R$ is the
hyperfinite $II_1$ factor with the canonical trace $\tau$.
\[
  [R: \Phi (R)]= {\rm rank} (A)^2, \quad H(\Phi)= \ln ({\rm rank}
  (A))
\]
\end {lemma}

The minimal projection of $A$ has the same trace.
Observe that $\otimes_{i=1}^{\infty} A$ contains a  hereditary
maximal abelian subalgebra of $R$, which is a stringent condition
for calculating the entropy.

We construct a variant of the above endomorphism in the next
section. The main technical tool in the construction is \cite {rP88}
R.Powers' binary shifts. We provide here the details of n-unitary
shifts generalized by \cite {mC87} M.Choda for the convenience of
the reader.

Let $n$ be a positive integer. We treat a pair of sets $Q$ and
$S$ of integers satisfying the following condition $(*)$ for
some integer $m$:
\[
(*)
\begin {cases}
Q= (i(1), i(2), \cdots, i(m)), &0 \leq i(1) < i(2) < \cdots
< i(m), \\
S= (j(1), j(2), \cdots, j(m)), & j(l)= 1, 2, \cdots, n-1,\\
&\text {\ for \ } l= 1, 2, \cdots, m.
\end {cases}
\]
\begin {definition}
A unital $*$-endomorphism $\Psi$ of $R$ is called an n-unitary
shift of $R$ if there is a unitary $u \in R$ satisfying the following:\\
(1)$u^n= 1$;\\
(2)$R$ is generated by $\{u, \Psi (u), \Psi^2 (u), \cdots, \}$;\\
(3)$\Psi^k (u) u= u \Psi^k (u)$ or $\Psi^k (u) u= \gamma u \Psi^k
(u)$
for all $k= 1,2,\cdots$, where $\gamma = \exp (2 \pi \sqrt {-1} / n)$.\\
(4) for each $(Q, S)$ satisfying $(*)$, there are an integer $k (\geq
0)$ and a nontrivial $\lambda \in \mathbb {T} = \{ \mu \in \mathbb {C};
| \mu | = 1 \}$ such that
\[
  \Psi^k (u) u(Q, S)= \lambda u(Q, S) \Psi^k (u),
\]
where  $u(Q, S)$ is defined by
\[
  u(Q, S) = {\Psi^{i(1)} (u)}^{j(1)} {\Psi^{i(2)} (u)}^{j(2)}
  \cdots  {\Psi^{i(m)} (u)}^{j(m)}.
\]
\end {definition}

The unitary $u$ is called a generator of $\Psi$. Put $S (\Psi; u)=
\{ k; \Psi^k (u) u= \gamma u \Psi^k (u) \}$. Note that the above
condition (2) gives some rigidity on $S (\Psi; u)$.
The Jones index $[R: \Psi (R)]$ is n.

One interesting example of $S_1= S (\Psi_1; u_1)$ is $\{
1,3,6,10,15, \cdots,$ $\frac {1} {2} l(l+1),  \cdots \}$, which
corresponds to the n-stream $\{ 0 1 0 1 0 0 1 0 0 0 1 0 0 0 0 1$ $ 0
0 0 0 0 1 \cdots \}$. It is pointed out that the relative commutant
$\Psi_1^k (R)' \cap R$ is always trivial for all $k$!

We define
\[
  S_2= \{ \frac{1}{2} l (l+1) \ | \ l = 1 \mod 3 \}, \quad
  S_3= \{ \frac{1}{2} l (l+1) \ | \ l = 2 \mod 3 \}
\]

\section {A Simplified Version} \label {S: simp}
\begin {theorem}
For any finite dimensional $C^*$-algebra $A$, there exists a $\Phi
\in End (R, \tau)$ such that the relative commutant $ \Phi^k (R)'
\cap R$ is isomorphic to $\otimes_{i=1}^k Z(A)$. $Z(A)$ is the
center of $A$. Here $R$ is the hyperfinite $II_1$ factor with the
trace $\tau$.
\end {theorem}

Since $A$ is finite dimensional, then $A$ can be decomposed as
a direct sum of finitely many matrix algebras,
\[
  A \simeq \oplus_{i=1}^j M_{a_i} (\mathbb {C}) \subseteq M_n
  (\mathbb {C}),
\]
where $n= \sum_i a_i$.
The minimal central projection of $A$ is not necessarily equivalent to
each other, though the minimal projection of $A$ is.
The trace vector on $Z (A)$ is 
\[ 
  [\frac{a_1}{n}, \ \frac{a_2}{n}, \ \cdots \ \frac{a_j}{n}]
\]

For each $i$, $M_{a_i} (\mathbb {C}) \subset A$ (not a unital
embedding) is generated by $p_i, q_i \in \mathcal {U}(\mathbb
{C}^{a_i})$ with:
\[
  p_i^{a_i}= q_i^{a_i}= 1_{M_{a_i}(\mathbb {C})}; \quad
  \gamma_i= \exp( 2 \pi \sqrt {-1} / a_i), \quad p_i q_i = \gamma_i  q_i p_i
\]
where $p_i= [ 1 \ \gamma_i \ \gamma_i^2 \cdots \gamma_i^{a_i-1} ]$
is the diagonal matrix in $M_{a_i} (\mathbb {C})$,
 and
$q_i$ is the permutation matrix in $M_{a_i} (\mathbb {C})$,
$( 1 \ 2 \ 3 \cdots a_i )$.

Define $v \in M_n( \mathbb {C})$ to be the permutation matrix:
\[
v= (a_1 \ \ (a_1+ a_2) \ \ (a_1+ a_2+ a_3) \ \ \cdots \ \
   (a_1+ a_2+ \cdots +a_j)).
\]
Then $v^j= 1$.
$A$ and $v$ do generate $M_n (\mathbb
{C})$. Therefore we can two describe $M_n (\mathbb {C})$
via $A$ and $v$ (or via $A$ and $r$ described below.)

Define $r:= s v s$, while
\[
  s= [0 \ 0 \ \cdots \ 0 \ 1_{a_1} \ 0 \ 0 \
  \cdots \ 0 \ 1_{a_1+ a_2} \cdots \ 0 \ 0 \cdots \ 0 \ 1_{a_1+ a_2+ \cdots+ a_j}]
  \in A.
\]
$A$ and $r$ generate $M_n (\mathbb {C})$.

\begin {lemma}
In fact, $<A, r> \simeq M_n (\mathbb {C})$ is of the form:
\[
  A+ A r A+ A r^2 A+ \cdots + A r^{j-1} A.
\]
\end {lemma}

\begin {proof}
It suffices to observe that
\begin{align*}
   &r= s v s= v s, \quad [s, v]= [s, r]=0\\
   &r^*= s v^*= v^* s, \quad r r^*= r^* r= s\\
   &r A r= r s A  r= r^2 r^* A r=  r^2 s (v^* A v) s \subset r^2 A\\
   &r^*= r^{j-1}, \quad r^j= s
\end {align*}
\end {proof}

On the other hand, define $w= \sum_{i=1}^j \gamma^{i-1} 1_{M_{a_i}
(\mathbb {C})}$, where 
\[
  \gamma= \exp (2 \pi \sqrt {-1} / j) \quad
   \gamma^j = 1
\]
Note that $w$ is in the center of $A$.
Two simple yet important observations are that:\\
(1) $Ad w$ acts trivially on $A$.\\
(2) $Ad w (r)= \gamma r$.

Now we  construct a tower of inclusion of finite dimensional $C^*$-
algebras $M_k$ with a trace $\tau$. The ascending union $M= \cup_{k
\in \mathbb {N}} M_k$ contains infinite copies of $M_n (\mathbb
{C})$, and thus of $A$. Number them respectively by
$r_1, A_1, w_1$, $r_2, A_2, w_2$, $r_3, A_3, w_3$, $\cdots$.
The key point in the construction is how the full matrix algebra 
$<r_k, A_k> \simeq M_n ({\mathbb C})$ is embedded in $M_k$.

We endow on this algebra the following properties:
\begin{align*}
  &[r_l, A_m]= 0, \ \text{if} \ l \not= m;\\
  &r_l r_m = \gamma r_l r_m,\  \text{if} \  |l- m| \in
  S_1= \{ 1, 3, 6, 10, 15, \cdots \},\\
  &r_l r_m = r_m r_l, \ \text{otherwise}.
\end{align*}
Here $\gamma= \exp (2 \pi \sqrt {-1} / j)$ as above.

Unlike in \cite {hH}, we add a twist in the relations between $A_l$
and $A_m$ when $|l-m| \in S_2 \cup S_3$.
\begin{align*}
  &p_{i, l} p_{i', m}= p_{i, m} p_{i', l};\\
   &q_{i, l} p_{i', m}= \gamma_i^{- \delta_{i, i'}} p_{i', m} q_{i, l}, \quad
  {\rm if \ } |l-m| \in S_2,\\
  &q_{i, l} p_{i', m}= p_{i', m} q_{i, l}, \quad
  {\rm if \ otherwise;} \\
  &q_{i, l} q_{i', m}= \gamma_i^{\delta_{i, i'}} q_{i', m} q_{i, l}, \quad
  {\rm if \ } |l-m| \in S_3\\
  &q_{i, l} q_{i', m}= q_{i', m} q_{i, l}, \quad
  {\rm if \ otherwise;}
\end{align*}
where $A_l= <p_{i, l}, q_{i, l}>_{i=1}^{j}$, and 
$A_m= <p_{i', m}, q_{i', m}>_{i'=1}^{j}$.

The construction is an induction process. We have handy the
embedding $A_1 \subseteq M_n (\mathbb {C})= M_1$, which is
isomorphic to the inclusion of $A \otimes 1_{M_n (\mathbb {C})}$
inside $M_n (\mathbb {C}) \otimes 1_{M_n (\mathbb {C})}$ equipped
with the trace $\frac{1}{n} Tr$.

Observe that $|2-1| = 1 \in S_2$. We would like to identify $A_2$ in
$\otimes^2 M_n (\mathbb {C})$ by a twist. $A_2$ is generated by
$p_{1, 2}, \cdots, p_{j, 2}$ and $q_{1, 2}, \cdots, q_{j, 2}$.

Put $p_{i, 2}= {\bf 1} \otimes p_i \in \otimes^2 M_n ({\mathbb C})$,
$q_{i, 2}= (q_{i} + {\bf 1}- 1_{M_{a_i} ({\mathbb C})}) \otimes q_i
\in \otimes^2 M_n ({\mathbb C})$. Note that $q_{i} + {\bf 1}-
1_{M_{a_i} ({\mathbb C})} \in {\mathcal U} ({\mathbb C}^n)$.  We
have:
\begin {align*}
    &p_{i, 2}^{a_i}= q_{i, 2}^{a_i}= {\bf 1} \otimes 1_{M_{a_i}
   ({\mathbb C})},\\
   &p_{i,2} q_{i,2}= \gamma_i q_{i,2} p_{i,2},\\
   &q_{i,2} p_{i,1}= \gamma_i^{-1} p_{i,1} q_{i, 2},\\
   &[p_{i,2}, p_{i,1}]= [q_{i,2}, q_{i,1}]= 0.
\end{align*}

$A_2$ is generated by $p_{1, 2}, \cdots, p_{j, 2}$ and $q_{1, 2},
\cdots, q_{j, 2}$.
\[
  A_2 \simeq A= \oplus_{i=1}^j M_{a_i} ({\mathbb C})
\]

Observe $|2-1|=1 \in S_1$. Define $r_2: = w \otimes r$. We have the
following properties:
\begin{align*}
   &[w, q_{i} + {\bf 1}- 1_{M_{a_i} ({\mathbb C})}]= 0\\
   &<A_2, r_2> \simeq M_n (\mathbb {C})\\
   &[A_1, r_2]= 0 \quad r_1 r_2= \gamma r_2 r_1\\
   &M_2:= <A_1, r_1, A_2, r_2> = \otimes^2 M_n ({\mathbb C})
\end{align*}
There is  a unique normalized trace $\tau$ on $M_2$.

Assume we have obtained $M_k= <A_1, r_1, A_2, r_2, \cdots, A_k,
r_k>$ equal to $ \otimes^k M_n ({\mathbb C})$ with the trace $\tau$.
We identify $M_k$ as $M_k \otimes 1_{M_n (\mathbb {C})}$ by sending
$x \in M_k$ to $x \otimes 1_{M_n (\mathbb {C})}$.

Define $A_{k+1}$ by its generators, $p_{i, k+1}$ and $q_{i, k+1}$,
$1 \leq i \leq j$:
\begin{align*}
   &p_{i, k+1} := {\bf 1} \otimes {\bf 1} \otimes \cdots \otimes {\bf 1} \otimes p_i\\
   &q_{i, k+1} := \\
   &[(q_{i} + {\bf 1}- 1_{M_{a_i} ({\mathbb C})})^{b_1}
    \otimes \cdots
    \otimes (q_{i} + {\bf 1}- 1_{M_{a_i} ({\mathbb C})})^{b_k}
    ] \cdot \\
   &[(p_{i} + {\bf 1}- 1_{M_{a_i} ({\mathbb C})})^{c_1}
    \otimes \cdots
    \otimes (p_{i} + {\bf 1}- 1_{M_{a_i} ({\mathbb C})})^{c_k}]
    \otimes q_i\\
    &b_l= 1, \quad {\rm if \ } |k+1-l| \in S_2; \quad
   b_l= 0, \quad {\rm otherwise}\\
   &c_l= 1, \quad {\rm if \ } |k+1-l| \in S_3; \quad
   c_l= 0, \quad {\rm otherwise}
\end{align*}

\begin{align*}
   &p_{i, k+1}^{a_i}= q_{i, k+1}^{a_i}= \otimes^k {\bf 1} \otimes 1_{M_{a_i}}({\mathbb
   C})\\
   &p_{i, k+1} q_{i, k+1}= \gamma_i q_{i, k+1} p_{i, k+1}
\end{align*}
Therefore $A_{k+1}$ is isomorphic to $A$.

The commutation relations is given below.
\begin{align*}
   &p_{i, k+1} p_{i', l}= p_{i', l} p_{i, k+1}\\
   &q_{i, k+1} p_{i', l}= \gamma_i^{-\delta_{i, i'}} p_{i', l} q_{i, k+1}, 
   \ {\rm if \ } |k+1-l| \in S_2\\
   &[q_{i, k+1}, p_{i', l}]= 0 , \ {\rm if \ } |k+1-l| \notin S_2\\
   &q_{i, k+1} q_{i', l}= \gamma_i^{\delta_{i, i'}} q_{i', l} q_{i, k+1}, 
   \ {\rm if \ } |k+1-l| \in S_3\\
   &[q_{i, k+1}, q_{i', l}]= 0 , \ {\rm if \ } |k+1-l| \notin S_3\\
   &A_{k+1} \cdot A_l= A_l \cdot A_{k+1}
\end{align*}

Define
\begin {align*}
   &r_{k+1}:= w^{d_1} \otimes w^{d_2} \otimes \cdots \otimes w^{d_k}
   \otimes r\\
   &d_l= 1, \quad {\rm if \ } |k+1-l| \in S_1; \quad
   d_l= 0, \quad {\rm otherwise}
\end{align*}
a

We have the following properties:
\begin{align*}
   &[w, q_{i} + {\bf 1}- 1_{M_{a_i} ({\mathbb C})}]= [w, p_{i} + 
{\bf 1}- 1_{M_{a_i} ({\mathbb C})}]= 0\\
   &<A_{k+1}, r_{k+1}> \simeq M_n (\mathbb {C})\\
   &[A_l, r_{k+1}]= 0 \quad 1 \leq l \leq k\\
    &r_{k+1}  r_l= \gamma r_l r_{k+1} \quad {\rm if \ } |k+1-l| \in S_1\\
    &r_{k+1}  r_l= r_l r_{k+1} \quad {\rm if \ } |k+1-l| \notin S_1\\
   &M_{k+1}:= <M_k, A_{k+1}, r_{k+1}> = \otimes^{k+1} M_n ({\mathbb C})
\end{align*}
There is  a unique normalized trace $\tau$ on $M_{k+1}$.

By induction we have constructed the ascending tower of finite
dimensional $C^*$-algebras with the desired properties.

We now explore some useful properties of the finite dimensional
$C^*$-algebra, $M_k$.

\begin {lemma}
For all $l$, $M_l$ is the linear span of the words, $x_1 \cdot x_2
\cdot x_3  \cdots  x_l$, where $x_k \in M_n (\mathbb {C})_k= <A_k,
r_k>$.
\end {lemma}

\begin {proof}
It suffices to prove $x_l \cdot x_k$ is in $M_k \cdot <A_l, r_l>= 
M_k \cdot M_n (\mathbb {C})_l$, where $k < l$.

\begin {align*}
&A_l A_k= A_k A_l \subset M_k A_l\\
&A_l r_l A_k= A_k A_l r_l \subset M_k A_l r_l\\
&p_{i, l} r_k= r_k p_{i, l} \in r_k A_l \subset M_k A_l
\end{align*}

Note that 
\begin{align*}
&[w, q_{i} + {\bf 1}- 1_{M_{a_i} ({\mathbb C})}]= [w, p_{i} + 
{\bf 1}- 1_{M_{a_i} ({\mathbb C})}]= 0;\\
&{\rm if \ } |l-k| \in S_2, {\rm \ then}\\
&q_{i, l} r_k= Ad ({\bf 1} \otimes {\bf 1} \otimes \cdots \otimes
(q_i + {\bf 1}- 1_{M_{a_i} ({\mathbb C})})) (r_k) \cdot
q_{i, l} \in M_k A_l;\\
&{\rm if \ } |l-k| \in S_3, {\rm \ then}\\
&q_{i, l} r_k= Ad ({\bf 1} \otimes {\bf 1} \otimes \cdots \otimes
(p_i + {\bf 1}- 1_{M_{a_i} ({\mathbb C})})) (r_k) \cdot q_{i, l} 
\in M_k A_l;\\
&{\rm if \ otherwise,} \quad q_{i, l} r_k = r_k q_{i, l} \in r_k A_l \subset M_k A_l.
\end{align*}

In short, 
\begin{align*}
&A_l r_k \subset M_k A_l;\\
&A_l A_k r_k= A_k A_l r_k \subset A_k M_k A_l \subset  M_k A_l;\\
&A_l r_l A_k r_k= A_k A_l r_l r_k= A_k A_l r_k r_l \subset A_k M_k
A_l r_l \subset M_k  A_l r_l.
\end {align*}
\end {proof}

\begin {lemma}
Consider the pair $(M, \tau)$ as described above and the GNS-
construction. Identify everything mentioned above as its image. We
$M''$ is the hyperfinite $II_1$ factor.
\end {lemma}

\begin{proof}
There is one and only one tracial state on $M_k$ for all $k \in
{\mathbb N}$. Hence the tracial state on $M$ is unique. Therefore
$M''$ is the hyperfinite $II_1$ factor, $R$.
\end {proof}

Define a unital *-endomorphism, $\Phi$, on $R$ to be the (right)
one-shift: i.e., sending $A_k$ to $A_{k+1}$, and sending $r_k$ to
$r_{k+1}$. We observe that $\Phi (R)$ is a $II_1$ factor and
\[
  [R: \Phi (R)]= n^2.
\]

\begin {lemma}
The relative commutant $\Phi^k (R)' \cap R$ is exactly
$\otimes_{i=1}^k Z(A)$, $Z(A)$ is the center of $A$.
\end {lemma}

\begin {proof}
Because of our decomposition in Lemma 3 and Lemma 4, $R$ can be
written as
\[
  (\sum_{i=0}^j A_1 r_1^{i} A_1) \cdot (\sum_{i=0}^j
  A_2 r_2^{i} A_2) \cdot \cdots \cdot (\sum_{i=0}^j A_k r_k^{i} A_k)
  \cdot \Phi^k (R).
\]

Assume $x \in R \cap \Phi^k (R)'$. $x$ can be written, as in Lemma 4,
of the following form:
\[
  x= \sum_{\vec {\alpha} \in \{ 0, \ 1,  \ \cdots, j-1 \}^k} y_1^{\vec {\alpha}} 
r_1^{c_1} z_1^{\vec
  {\alpha}} y_2^{\vec {\alpha}} r_2^{c_2} z_2^{\vec {\alpha}} 
\cdots y_k^{\vec {\alpha}} r_1^{c_k}
  z_k^{\vec {\alpha}}
  \cdot y^{\vec {\alpha}},
\]
where $\vec {\alpha}= (c_1, c_2, \dots, c_k)$ is a multi-index and
$y^{\vec {\alpha}}$ is in $\Phi^k (R)$. Note that $\Phi^k (R)$ is
the weak closure of $\{ \Phi^k (M_i) \}_{i=1}^{\infty}$.

For every $\epsilon > 0$,  there exists an integer $i \in {\mathbb
N}$ such that
\begin{align*}
   &\forall \ {\vec {\alpha}}, \ z^{\vec {\alpha}} \in \Phi^k (M_i)
   \subset <A_{k+1}, r_{k+1}, \cdots, A_{k+i}, r_{k+i}>\\
   &\| x-\sum_{\vec {\alpha} \in \{ 0, \ 1,  \ \cdots, j-1 \}^k} y_1^{\vec {\alpha}} 
r_1^{c_1} z_1^{\vec
  {\alpha}}
  y_2^{\vec {\alpha}} r_2^{c_2} z_2^{\vec {\alpha}} \cdots y_k^{\vec {\alpha}} r_1^{c_k}
  z_k^{\vec {\alpha}}
  \cdot z^{\vec {\alpha}} \|_{2, \tau} < \delta\\
  &\delta= (\sqrt{\frac{j}{n}})^k \epsilon
\end{align*}

Put $L= l( l+1)/2+1$ for some integer $l > k+1$ and $l= 0 \mod 3$.
We have the following properties:
\begin {align*}
&[r_L, A_1]= [r_L, A_2]= \cdots = [r_L, A_{k+i}]= 0\\
&[r_L, r_2]= [r_L, r_3]= \cdots = [r_L, r_{k+i}]= 0\\
&r_L r_1 = \gamma r_1 r_L, \quad r_L r_1^{c_1}
   = \gamma^{c_1} r_1 r_L\\
&r_L r_L^*= r_L^* r_L= s_L\\
&r_L r_1 r_L^*= \gamma r_1 s_L, \quad r_L r_1^{c_1} r_L^*= \gamma^{c_1} r_1 s_L\\
&{\rm for \ } 0 \leq m \leq {j-1}, \quad r_L^m r_1^{c_1} (r_L^*)^m= \gamma^{c_1 m} r_1 s_l\\
&[s_L, A_1]= [s_L, A_2]= \cdots = [s_L, A_{k+i}]= 0\\
&[s_L, r_1]= [s_L, r_2]= \cdots = [s_L, r_{k+i}]= 0
\end {align*}

Therefore we claim:
\begin {align*}
&\| (x-\sum_{\vec {\alpha} \in \{ 0, \ 1,  \ \cdots, j-1 \}^k}
y_1^{\vec {\alpha}} r_1^{c_1} z_1^{\vec
  {\alpha}}
  y_2^{\vec {\alpha}} r_2^{c_2} z_2^{\vec {\alpha}} \cdots y_k^{\vec {\alpha}} r_1^{c_k}
  z_k^{\vec {\alpha}}
  \cdot z^{\vec {\alpha}}) s_L \|_{2, \tau}=\\
&\| (x-\sum_{\vec {\alpha} \in \{ 0, \ 1,  \ \cdots, j-1 \}^k}
y_1^{\vec {\alpha}} r_1^{c_1} z_1^{\vec
  {\alpha}}
 \cdots y_k^{\vec {\alpha}} r_1^{c_k}
  z_k^{\vec {\alpha}}
  \cdot z^{\vec {\alpha}}) \frac{1}{j} \sum_{m=0}^{j-1} r_L^m (r_L^*)^m \|_{2, \tau}=\\
&\frac{1} {j} \| \sum_{\vec {\alpha}} \sum_m (r_L^m x {r_L^*}^m-
y_1^{\vec {\alpha}} r_L^m r_1^{c_1} (r_L^*)^m z_1^{\vec
  {\alpha}}
  y_2^{\vec {\alpha}} r_2^{c_2} z_2^{\vec {\alpha}} \cdots y_k^{\vec {\alpha}} r_1^{c_k}
  z_k^{\vec {\alpha}}
  \cdot z^{\vec {\alpha}}) \|_{2, \tau}=\\
&\frac{1} {j} \| \sum_{\vec {\alpha}} \sum_m (x - y_1^{\vec
{\alpha}}  r_1^{c_1 m} z_1^{\vec
  {\alpha}}
  y_2^{\vec {\alpha}} r_2^{c_2} z_2^{\vec {\alpha}} \cdots y_k^{\vec {\alpha}} r_1^{c_k}
  z_k^{\vec {\alpha}}
  \cdot z^{\vec {\alpha}}) s_L \|_{2, \tau}=\\
&\| (x-\sum_{\vec {\alpha} \in \{ 0, \ 1, \cdots, j-1 \}^k, c_1= 0}
y_1^{\vec {\alpha}} z_1^{\vec
  {\alpha}}
  y_2^{\vec {\alpha}} r_2^{c_2} z_2^{\vec {\alpha}} \cdots y_k^{\vec {\alpha}} r_1^{c_k}
  z_k^{\vec {\alpha}}
  \cdot z^{\vec {\alpha}}) s_L \|_{2, \tau}=\\
&\sqrt{\frac{j}{n}} \| x-\sum_{\vec {\alpha} \in \{ 0, \ 1, \cdots,
j-1 \}^k, c_1= 0} y_1^{\vec {\alpha}} z_1^{\vec
  {\alpha}}
  y_2^{\vec {\alpha}} r_2^{c_2} z_2^{\vec {\alpha}} \cdots y_k^{\vec {\alpha}} r_1^{c_k}
  z_k^{\vec {\alpha}}
  \cdot z^{\vec {\alpha}}  \|_{2, \tau}\\
&{\rm since \ } \{ x, \ y_1^{\vec {\alpha}}, \ z_1^{\vec
  {\alpha}}, \ y_2^{\vec {\alpha}}, \ r_2^{c_2}, \ z_2^{\vec {\alpha}}, \cdots, y_k^{\vec{\alpha}},
  \ r_1^{c_k}, \ z_k^{\vec {\alpha}}, \ z^{\vec {\alpha}} \} \subset \{ s_L, \ r_L, \ A_L
  \}'\\
&{\rm and \ } \tau (s_L)= \frac{j}{n}.
\end {align*}
Note that $\{ s_L, \ r_L, \ A_L \}'' = M_n({\mathbb C})_L$ is a
type $I$ factor \cite{sP83}.

By induction,
\begin{align*}
&\| x-\sum_{\vec {\alpha} \in \{ 0, \ 1, \cdots, j-1 \}^k, \ c_1= 0}
y_1^{\vec {\alpha}} z_1^{\vec
  {\alpha}}
  y_2^{\vec {\alpha}} r_2^{c_2} z_2^{\vec {\alpha}} \cdots y_k^{\vec {\alpha}} r_1^{c_k}
  z_k^{\vec {\alpha}}
  \cdot z^{\vec {\alpha}}  \|_{2, \tau} < \sqrt{\frac{n}{j}} \delta\\
&\| x-\sum_{\vec {\alpha} \in \{ 0, \ 1, \cdots, j-1 \}^k, \ c_1=
c_2= 0} y_1^{\vec {\alpha}} z_1^{\vec
  {\alpha}}
  y_2^{\vec {\alpha}} z_2^{\vec {\alpha}} \cdots y_k^{\vec {\alpha}} r_1^{c_k}
  z_k^{\vec {\alpha}}
  \cdot z^{\vec {\alpha}}  \|_{2, \tau} < (\sqrt{\frac{n}{j}})^2 \delta\\
&\cdots \\
&\| x-\sum_{\vec {\alpha} \in \{ 0 \}^k} y_1^{\vec {\alpha}}
z_1^{\vec
  {\alpha}}
  y_2^{\vec {\alpha}} z_2^{\vec {\alpha}} \cdots y_k^{\vec {\alpha}}
  z_k^{\vec {\alpha}}
  \cdot z^{\vec {\alpha}}  \|_{2, \tau} < (\sqrt{\frac{n}{j}})^k \delta= \epsilon\\
\end{align*}

Put $L_1= l_1( l_1+1)/2+1$ for some integer $l_1 > k+1$ and $l_1= 1
\mod 3$. We have the following properties:
\begin{align*}
   &U_{L_1}:= \sum_{m_1=0}^{j-1} q_{m_1, L_1}\\
   &U_{L_1}^n= {\bf 1}\\
   &[U_{L_1}, \Phi (M_{k+i-1})]= 0
\end{align*}

Similarly, put $L_2= l_2( l_2+1)/2+1$ for some integer $l_2 > k+1$
and $l_2= 2 \mod 3$. We have the following properties:
\begin{align*}
   &U_{L_2}:= \sum_{m_2=0}^{j-1} q_{m_2, L_2}\\
   &U_{L_2}^n= {\bf 1}\\
   &[U_{L_2}, \Phi (M_{k+i-1})]= 0
\end{align*}

Therefore
\begin{align*}
&\| x-\sum_{\vec {\alpha} \in \{ 0 \}^k} y_1^{\vec {\alpha}}
z_1^{\vec
  {\alpha}}
  y_2^{\vec {\alpha}} z_2^{\vec {\alpha}} \cdots y_k^{\vec {\alpha}}
  z_k^{\vec {\alpha}}
  \cdot z^{\vec {\alpha}}  \|_{2, \tau} =\\
&\| x-\sum_{\vec {\alpha} \in \{ 0 \}^k} y_1^{\vec {\alpha}}
z_1^{\vec
  {\alpha}}
  y_2^{\vec {\alpha}} z_2^{\vec {\alpha}} \cdots y_k^{\vec {\alpha}}
  z_k^{\vec {\alpha}}
  \cdot z^{\vec {\alpha}}  \frac {1}{n} \sum_{m_1=0}^{n-1} U_{L_1}^{m_1} {U_{L_1}^*}^{m_1} \|_{2,
  \tau}=\\
&\| x-\sum_{\vec {\alpha} \in \{ 0 \}^k} \frac {1}{n}
\sum_{m_1=0}^{n-1} U_{L_1}^{m_1} (y_1^{\vec {\alpha}} z_1^{\vec
  {\alpha}}) {U_{L_1}^*}^{m_1}
  y_2^{\vec {\alpha}} z_2^{\vec {\alpha}} \cdots y_k^{\vec {\alpha}}
  z_k^{\vec {\alpha}}
  \cdot z^{\vec {\alpha}}  \|_{2, \tau}=\\
&\| x-\sum_{\vec {\alpha} \in \{ 0 \}^k} \frac {1}{n^2}
\sum_{m_1, m_2=0}^{n-1} U_{L_1}^{m_1} (y_1^{\vec {\alpha}} z_1^{\vec
  {\alpha}}) {U_{L_1}^*}^{m_1}
\cdots y_k^{\vec {\alpha}}
  z_k^{\vec {\alpha}}
  \cdot z^{\vec {\alpha}}  
 U_{L_2}^{m_2} {U_{L_2}^*}^{m_2}
  \|_{2, \tau}=\\
&\| x-\sum_{\vec {\alpha} \in \{ 0 \}^k} \frac {1}{n^2} \sum_{m_1,
m_2=0}^{n-1} U_{L_2}^{m_2} U_{L_1}^{m_1} (y_1^{\vec {\alpha}}
z_1^{\vec
  {\alpha}}) {U_{L_1}^*}^{m_1} {U_{L_2}^*}^{m_2}
\cdots y_k^{\vec {\alpha}}
  z_k^{\vec {\alpha}}
  \cdot z^{\vec {\alpha}}
  \|_{2, \tau}= \\
&\| x-\sum_{\vec {\alpha} \in \{ 0 \}^k} x_1^{\vec {\alpha}}
  y_2^{\vec {\alpha}} z_2^{\vec {\alpha}} \cdots y_k^{\vec {\alpha}}
  z_k^{\vec {\alpha}}
  \cdot z^{\vec {\alpha}}
  \|_{2, \tau}
\end{align*}

Observe that
\[
  x_1^{\vec {\alpha}}:= \frac {1}{n^2} \sum_{m_1, m_2=0}^{n-1} U_{L_2}^{m_2}
 U_{L_1}^{m_1} (y_1^{\vec {\alpha}} z_1^{\vec
  {\alpha}}) {U_{L_1}^*}^{m_1} {U_{L_2}^*}^{m_2}
\]
is the conditional expectation of $y_1^{\vec {\alpha}} z_1^{\vec
  {\alpha}}$ on $Z(A_1)$, the center of $A_1$.

By induction,
\begin{align*}
&\| x-\sum_{\vec {\alpha} \in \{ 0 \}^k} x_1^{\vec {\alpha}}
  (y_2^{\vec {\alpha}} z_2^{\vec {\alpha}}) \cdots y_k^{\vec {\alpha}}
  z_k^{\vec {\alpha}}
  \cdot z^{\vec {\alpha}}
  \|_{2, \tau}=\\
&\| x-\sum_{\vec {\alpha} \in \{ 0 \}^k} x_1^{\vec {\alpha}}
  x_2^{\vec {\alpha}}  (y_3^{\vec {\alpha}}
  z_3^{\vec {\alpha}}) \cdots y_k^{\vec {\alpha}}
  z_k^{\vec {\alpha}}
  \cdot z^{\vec {\alpha}}
  \|_{2, \tau}=\\
&\cdots\\
&\| x-\sum_{\vec {\alpha} \in \{ 0 \}^k} x_1^{\vec {\alpha}}
  x_2^{\vec {\alpha}}  \cdots x_k^{\vec {\alpha}}
  \cdot z^{\vec {\alpha}}
  \|_{2, \tau} < \epsilon
\end{align*}
where $x_1^{\vec {\alpha}} \in Z(A_1)$, $x_2^{\vec {\alpha}} \in
Z(A_2)$, $\cdots$, $x_k^{\vec {\alpha}} \in Z(A_k)$.

Note that the von Neumann algebra  $\{ x, \ x_1^{\vec {\alpha}}, \
x_2^{\vec {\alpha}}, \ x_3^{\vec {\alpha}}, \cdots x_k^{\vec
{\alpha}} \}''$ commutes with $\Phi^k (M)$, which is a $II_1$
factor. Any element in the former von Neumann algebra has a scalar
conditional expectation onto $\Phi^k (M)$. In short, the former von
Neumann algebra and $\Phi^k (M)$ are mutually orthogonal.

According to \cite {sP83}, we have:
\begin{align*}
  &\| x-\sum_{\vec {\alpha} \in \{ 0 \}^k} x_1^{\vec {\alpha}}
  x_2^{\vec {\alpha}} x_3^{\vec {\alpha}} \cdots x_k^{\vec {\alpha}}
  \cdot \tau (z^{\vec {\alpha}})  \|_{2, \tau} < \epsilon\\
  &\sum_{\vec {\alpha} \in \{ 0 \}^k} x_1^{\vec {\alpha}}
  x_2^{\vec {\alpha}} x_3^{\vec {\alpha}} \cdots x_k^{\vec {\alpha}}
  \cdot \tau (z^{\vec {\alpha}})  \in Z(A_1) \cdot Z(A_2) \cdots Z(A_k)= \otimes^k
  Z(A)
\end{align*}

\end {proof}

\begin{corollary}
The entropy of the endomorphism $\Phi$  with the domain 
restricted on
$\otimes_{i=1}^{\infty} Z(A)$ is equivalent to
\[
  \sum_{i=1}^j -\frac{a_i}{n} \ln (\frac{a_i}{n}) \leq \ln n
\]
\end{corollary}
Take $A= M_2 ({\mathbb C}) \oplus M_2 ({\mathbb C}) \subset M_4
({\mathbb C})$. The associated endomorphism has index equal to $4^2=
16$. Yet the entropy restricted on the relative commutant algebra
gives $\ln 2$.

\section {Main Theorem} \label {S:main}

\begin{theorem}
For any finite dimensional $C^*$-algebra $A$ with a trace vector
$\vec s$ whose entries are rational numbers, we give an endomorphism
$\Phi$ of the hyperfinite $II_1$ factor $R$  such that:
\[
  \forall \ k \ \in \ \mathbb {N}, \
  \Phi^k (R)' \cap R= \otimes^k A.
\]
The canonical trace $\tau$ on $R$ extends the trace vector $\vec s$
on $A$.
\end{theorem}

The proof is nothing but a tedious generalization of Theorem 1.

$A$ is characterized by its trace vector $\vec s$ and its dimension
vector $\vec t$.
\begin {align*}
   &\vec s= [\frac{e_1}{f_1}, \ \frac{e_2}{f_2}, \cdots,
   \frac{e_j}{f_j}]\\
   &\vec t= [a_1, \ a_2, \cdots, a_j]\\
   &\frac{e_1}{f_1} a_1+ \frac{e_2}{f_2} a_2+ \cdots + \frac{e_j}{f_j} a_j=
   1\\
   &e_1,\  f_1, \  a_1, \  e_2,\  f_2, \  a_2, \cdots, e_j, \ f_j, \ a_j \in
   {\mathbb N}
\end {align*}

Put $n= f_1 f_2 \cdots f_j$. We can embed $A$ into $B \subset M_n
({\mathbb C})$ via
\begin{align*}
  &A \simeq
  M_{a_1} ({\mathbb C}) \oplus M_{a_2} ({\mathbb C})
  \oplus \cdots \oplus M_{a_j} ({\mathbb C}) \subset\\
  &M_{a_1} ({\mathbb C}) \otimes M_{a_1'} ({\mathbb C})
   \oplus M_{a_2} ({\mathbb C}) \otimes M_{a_2'} ({\mathbb C})
  \oplus \cdots \oplus M_{a_j} ({\mathbb C}) \otimes M_{a_j'} ({\mathbb
  C})\\
  &= B \subset M_n ({\mathbb C})\\
  &{\rm where \ } a_1'= \frac{n
  e_1}{f_1}, \ a_2'= \frac{n e_2}{f_2}, \cdots, a_j'=
  \frac{n e_j}{f_j}
\end{align*}

For each $i$, $M_{a_i} (\mathbb {C}) \subset A \subset B$ (the
former being not a unital embedding) is generated by $p_i, q_i \in
\mathcal {U}(\mathbb {C}^{a_i})$ with:
\[
  p_i^{a_i}= q_i^{a_i}= 1_{M_{a_i}(\mathbb {C})}; \quad
  \gamma_i= \exp( 2 \pi \sqrt {-1} / a_i), \quad p_i q_i = \gamma_i  q_i p_i
\]
where $p_i= [ 1 \ \gamma_i \ \gamma_i^2 \cdots \gamma_i^{a_i-1} ]$
is the diagonal matrix in $M_{a_i} (\mathbb {C})$,
 and
$q_i$ is the permutation matrix in $M_{a_i} (\mathbb {C})$, $( 1 \ 2
\ 3 \cdots a_i )$.

For each $i$, $M_{a_i'} (\mathbb {C}) \subset A' \cap B \subset B$
(the former being not a unital embedding) is generated by $p_i',
q_i' \in \mathcal {U}(\mathbb {C}^{a_i'})$ with:
\[
  p_i'^{a_i'}= q_i'^{a_i'}= 1_{M_{a_i'}(\mathbb {C})}; \quad
  \gamma_i'= \exp( 2 \pi \sqrt {-1} / a_i'), \quad p_i' q_i' = \gamma_i'  q_i'
  p_i'
\]
where $p_i'= [ 1 \ \gamma_i' \ \gamma_i'^2 \cdots \gamma_i'^{a_i-1}
]$ is the diagonal matrix in $M_{a_i'} (\mathbb {C})$,
 and
$q_i'$ is the permutation matrix in $M_{a_i'} (\mathbb {C})$, $( 1 \
2 \ 3 \cdots a_i' )$.

Define $v \in M_n( \mathbb {C})$ to be the permutation matrix:
\[
v= (a_1 a_1' \  \ (a_1 a_1'+ a_2 a_2') \ \
\cdots \ \  (a_1 a_1'+ a_2 a_2'+ \cdots +a_j a_j'))
\]
Then $v^j= 1$. $B$ and $v$  generate $M_n (\mathbb {C})$.

Define $r:= s v s$, while
\[
  s= [0 \ 0 \ \cdots \ 0 \ 1_{a_1 a_1'} \ 0 \ 0 \
  \cdots \ 0 \ 1_{a_1 a_1'+ a_2 a_2'} \cdots \ 0 \ 0 \cdots \ 0 \ 1_{a_1 a_1'+ a_2 a_2'+
  \cdots+ a_j a_j'}]
\]
is a diagonal matrix in $M_n ({\mathbb C})$.  Thus $B$ and $r$
generates $M_n (\mathbb {C})$.

\begin {lemma}
In fact, $<B, r> \simeq M_n (\mathbb {C})$ is of the form:
\[
  B+ B r B+ B r^2 B+ \cdots + B r^{j-1} B.
\]
\end {lemma}

\begin {proof}
It suffices to observe that
\begin{align*}
   &r= s v s= v s, \quad [s, v]= [s, r]=0\\
   &r^*= s v^*= v^* s, \quad r r^*= r^* r= s\\
   &r B r= r s B  r= r^2 r^* B r=  r^2 s (v^* B v) s \subset r^2 B\\
   &r^*= r^{j-1}, \quad r^j= s
\end {align*}
\end {proof}

On the other hand, define $w= \sum_{i=1}^j \gamma^{i-1} 1_{M_{a_i}
(\mathbb {C})}$, where 
\[
  \gamma= \exp (2 \pi \sqrt {-1} / j) \quad
  \gamma^j = 1. 
\]
Note that $w$ is in the center of $B$.
Two simple yet important observations are that:\\
(1) $Ad w$ acts trivially on $B$.\\
(2) $Ad w (r)= \gamma r$.

Now we  construct a tower of inclusion of finite dimensional $C^*$-
algebras $M_k$ with a trace $\tau$. The ascending union $M= \cup_{k
\in \mathbb {N}} M_k$ contains infinite copies of $M_n (\mathbb
{C})$, and thus of $B$. Number them respectively by $r_1, B_1, w_1$,
$r_2, B_2, w_2$, $r_3, B_3, w_3$, $\cdots$.

We endow on this algebra the following properties:
\begin{align*}
  &[r_l, B_m]= 0, \ \text{if} \ l \not= m;\\
  &r_l r_m = \gamma r_l r_m,\  \text{if} \  |l- m| \in
  S_1= \{ 1, 3, 6, 10, 15, \cdots \},\\
  &r_l r_m = r_m r_l, \ \text{otherwise}.
\end{align*}
Here $\gamma= \exp (2 \pi \sqrt {-1} / j)$ as above.

Unlike in \cite {hH}, we add a twist in the relations between $A_l'$
and $A_m'$ when $|l-m| \in S_2 \cup S_3$.
\begin{align*}
  &p_{i, l}' p_{i_0, m}'= p_{i_0, m}' p_{i, l}';\\
   &q_{i, l}' p_{i_0, m}'= \gamma_i'^{-\delta_{i, i_0}} p_{i_0, m}' q_{i, l}', \quad
  {\rm if \ } |l-m| \in S_2,\\
  &q_{i, l}' p_{i_0, m}'= p_{i_0, m}' q_{i, l}', \quad
  {\rm if \ otherwise;} \\
  &q_{i, l}' q_{i_0, m}'= \gamma_i'^{\delta_{i, i_0}} q_{i_0, m}' q_{i, l}', \quad
  {\rm if \ } |l-m| \in S_3\\
  &q_{i, l}' q_{i_0, m}'= q_{i_0, m}' q_{i, l}', \quad
  {\rm if \ otherwise;}
\end{align*}
where $A_l'= <p_{i, l}', q_{i, l}'>_{i=1}^{j}$, and $A_m'= <p_{i_0,
m}', q_{i_0, m}'>_{i_0=1}^{j}$.

Unlike in the above section, there is no twist in the relations
between $A_l$ and $A_m$.
\begin{align*}
  &p_{i, l} p_{i_0, m}= p_{i_0, m} p_{i, l};\\
   &p_{i, l} q_{i_0, m}= q_{i_0, m} p_{i, l};\\
  &q_{i, l} q_{i_0, m}= q_{i_0, m} q_{i, l}
\end{align*}
where $A_l= <p_{i, l}, q_{i, l}>_{i=1}^{j}$, and $A_m= <p_{i_0, m},
q_{i_0, m}>_{i_0=1}^{j}$.

The construction is an induction process. We have handy the
embedding $B_1 \subseteq M_n (\mathbb {C})= M_1$, which is
isomorphic to the inclusion of $B \otimes 1_{M_n (\mathbb {C})}$
inside $M_n (\mathbb {C}) \otimes 1_{M_n (\mathbb {C})}$ equipped
with the trace $\frac{1}{n} Tr$.

Observe that $|2-1| = 1 \in S_2$. We would like to identify $B_2$ in
$\otimes^2 M_n (\mathbb {C})$ by a twist. $B_2$ is generated by
$p_{1, 2}, \cdots, p_{j, 2}$,  $p_{1, 2}', \cdots, p_{j, 2}'$,
$q_{1, 2}, \cdots, q_{j, 2}$,and $q_{1, 2}', \cdots, q_{j, 2}'$.

Put $p_{i, 2}= {\bf 1} \otimes p_i \in \otimes^2 M_n ({\mathbb C})$,
$q_{i, 2}= {\bf 1} \otimes q_i \in \otimes^2 M_n ({\mathbb C})$. Put
$p_{i, 2}'= {\bf 1} \otimes p_i' \in \otimes^2 M_n ({\mathbb C})$,
$q_{i, 2}'= (q_i' + {\bf 1}- 1_{M_{a_i'} ({\mathbb C})}) \otimes
q_i' \in \otimes^2 M_n ({\mathbb C})$. Note that $q_i' + {\bf 1}-
1_{M_{a_i'} ({\mathbb C})} \in {\mathcal U} ({\mathbb C}^n)$.  We
have:
\begin {align*}
    &p_{i, 2}'^{a_i'}= q_{i, 2}'^{a_i'}= {\bf 1} \otimes 1_{M_{a_i'}
   ({\mathbb C})},\\
   &p_{i,2}' q_{i,2}'= \gamma_i' q_{i,2}' p_{i,2}',\\
   &q_{i,2}' p_{i,1}'= \gamma_i'^{-1} p_{i,1}' q_{i, 2}',\\
   &[p_{i,2}', p_{i,1}']= [q_{i,2}', q_{i,1}']= 0.
\end{align*}

$B_2$ is generated by 
\[
  p_{1, 2}, \cdots, p_{j, 2},  \ p_{1, 2}',
\cdots, p_{j, 2}', \ q_{1, 2}, \cdots, q_{j, 2}, \ q_{1, 2}',
\cdots, q_{j, 2}'
\]
Thus
\[
  B_2 \simeq B= \oplus_{i=1}^j M_{a_i a_i'} ({\mathbb C}) \subset M_n
  ({\mathbb C})
\]

Observe $|2-1|=1 \in S_1$. Define $r_2: = w \otimes r$. We have the
following properties:
\begin{align*}
   &[w, q_i' + {\bf 1}- 1_{M_{a_i'} ({\mathbb C})}]= 0\\
   &<B_2, r_2> \simeq M_n (\mathbb {C})\\
   &[B_1, r_2]= 0 \quad r_1 r_2= \gamma r_2 r_1\\
   &M_2= <B_1, r_1, B_2, r_2> = \otimes^2 M_n ({\mathbb C})
\end{align*}
There is  a unique normalized trace $\tau$ on $M_2$.

Assume we have obtained $M_k= <B_1, r_1, B_2, r_2, \cdots, B_k,
r_k>$ equal to $ \otimes^k M_n ({\mathbb C})$ with the trace $\tau$.
We identify $M_k$ as $M_k \otimes 1_{M_n (\mathbb {C})}$ by sending
$x \in M_k$ to $x \otimes 1_{M_n (\mathbb {C})}$.

Define $B_{k+1}$ by its generators: $p_{i, k+1}$, $p_{i, k+1}'$,
$q_{i, k+1}$, and $q_{i, k+1}'$, for all $1 \leq i \leq j$:
\begin{align*}
   &p_{i, k+1} := {\bf 1} \otimes {\bf 1} \otimes \cdots \otimes {\bf 1} \otimes p_i\\
   &q_{i, k+1} := {\bf 1} \otimes {\bf 1} \otimes \cdots \otimes {\bf 1} \otimes q_i\\
   &p_{i, k+1}' := {\bf 1} \otimes {\bf 1} \otimes \cdots \otimes {\bf 1} \otimes p_i'\\
   &q_{i, k+1}' := \\
   &[(q_i' + {\bf 1}- 1_{M_{a_i'} ({\mathbb C})})^{b_1}
    \otimes \cdots
    \otimes (q_i' + {\bf 1}- 1_{M_{a_i'} ({\mathbb C})})^{b_k}
    ] \cdot \\
   &[(p_i' + {\bf 1}- 1_{M_{a_i'} ({\mathbb C})})^{c_1}
    \otimes \cdots
    \otimes (p_i' + {\bf 1}- 1_{M_{a_i'} ({\mathbb C})})^{c_k}]
    \otimes q_i'\\
    &{\rm For \ all \ } 1 \leq l \leq k:\\
    &b_l= 1, \quad {\rm if \ } |k+1-l| \in S_2; \quad
   b_l= 0, \quad {\rm otherwise}\\
   &c_l= 1, \quad {\rm if \ } |k+1-l| \in S_3; \quad
   c_l= 0, \quad {\rm otherwise}
\end{align*}

We have:
\begin{align*}
   &p_{i, k+1}^{a_i}= q_{i, k+1}^{a_i}= \otimes^k {\bf 1} \otimes 1_{M_{a_i}}({\mathbb
   C})\\
   &p_{i, k+1} q_{i, k+1}= \gamma_i q_{i, k+1} p_{i, k+1}\\
    &p_{i, k+1}'^{a_i'}= q_{i, k+1}'^{a_i'}= \otimes^k {\bf 1} \otimes 1_{M_{a_i'}}({\mathbb
   C})\\
   &p_{i, k+1}' q_{i, k+1}'= \gamma_i' q_{i, k+1}' p_{i, k+1}'
\end{align*}
Therefore $B_{k+1}$ is isomorphic to $B$.

The commutation relations is given below.
\begin{align*}
   &p_{i, k+1} p_{i_0, l}= p_{i_0, l} p_{i, k+1}\\
   &[p_{i, k+1}, q_{i_0, l}]= 0 , \ {\rm for \ all \ }  l <k+1\\
   &[q_{i, k+1}, q_{i_0, l}]= 0 \\
   &A_{k+1} \cdot A_l= A_l \cdot A_{k+1}\\
   &p_{i, k+1}' p_{i_0, l}'= p_{i_0, l}' p_{i, k+1}'\\
   &q_{i, k+1}' p_{i_0, l}'= \gamma_i'^{-\delta_{i, i_0}} 
   p_{i_0, l}' q_{i, k+1}', \ {\rm if \ } |k+1-l| \in S_2\\
   &[q_{i, k+1}', p_{i_0, l}']= 0 , \ {\rm if \ } |k+1-l| \notin S_2\\
   &q_{i, k+1}' q_{i_0, l}'= \gamma_i'^{\delta_{i, i_0}} q_{i_0, l}' q_{i, k+1}', 
   \ {\rm if \ } |k+1-l| \in S_3\\
   &[q_{i, k+1}', q_{i_0, l}']= 0 , \ {\rm if \ } |k+1-l| \notin S_3\\
   &A_{k+1}' \cdot A_l'= A_l' \cdot A_{k+1}'
\end{align*}

Define
\begin {align*}
   &r_{k+1}:= w^{d_1} \otimes w^{d_2} \otimes \cdots \otimes w^{d_k}
   \otimes r\\
   &d_l= 1, \quad {\rm if \ } |k+1-l| \in S_1; \quad
   d_l= 0, \quad {\rm otherwise}
\end{align*}

We have the following properties:
\begin{align*}
   &[w, q_i' + {\bf 1}- 1_{M_{a_i'} ({\mathbb C})}]= [w, p_i' + {\bf 1}- 1_{M_{a_i'} ({\mathbb C})}]= 0\\
   &<B_{k+1}, r_{k+1}> \simeq M_n (\mathbb {C})\\
   &[B_l, r_{k+1}]= 0 \quad 1 \leq l \leq k\\
    &r_{k+1}  r_l= \gamma r_l r_{k+1} \quad {\rm if \ } |k+1-l| \in S_1\\
    &r_{k+1}  r_l= r_l r_{k+1} \quad {\rm if \ } |k+1-l| \notin S_1\\
   &M_{k+1}= <M_k, B_{k+1}, r_{k+1}> = \otimes^{k+1} M_n ({\mathbb C})
\end{align*}
There is  a unique normalized trace $\tau$ on $M_{k+1}$.

By induction we have constructed the ascending tower of finite
dimensional $C^*$-algebras with the desired properties.

We now explore some useful properties of the finite dimensional
$C^*$-algebra, $M_k$.

\begin {lemma}
For all $k$, $M_k$ is the linear span of the words, $x_1 \cdot x_2
\cdot x_3  \cdots  x_k$, where $x_j \in M_n (\mathbb {C})_j= <B_j,
r_j>$.
\end {lemma}

\begin {proof}
It suffices to prove $x_l \cdot x_j$ is in $M_j \cdot <B_l, r_l>= M_j \cdot 
M_n (\mathbb {C})_l$, where $j < l$.

\begin {align*}
&B_l B_j= B_j B_l \subset M_j B_l;\\
&B_l r_l B_j= B_j B_l r_l \subset M_j B_l r_l;\\
&p_{i, l} r_j= r_j p_{i, l} \in r_j B_l \subset M_j B_l;\\
&q_{i, l} r_j= r_j q_{i, l} \in r_j B_l \subset M_j B_l;\\
&p_{i, l}' r_j= r_j p_{i, l}' \in r_j B_l \subset M_j B_l;\\
&{\rm if \ } |l-j| \in S_2, {\rm \ then}\\
&q_{i, l}' r_j= Ad ({\bf 1} \otimes {\bf 1} \otimes \cdots \otimes
(q_i' + {\bf 1}- 1_{M_{a_i'} ({\mathbb C})})) (r_j) \cdot
q_{i, l}' \in M_j B_l;\\
&{\rm if \ } |l-j| \in S_3, {\rm \ then}\\
&q_{i, l}' r_j= Ad ({\bf 1} \otimes {\bf 1} \otimes \cdots \otimes
(p_i' + {\bf 1}- 1_{M_{a_i'} ({\mathbb C})})) (r_j) \cdot q_{i, l}' \in M_j B_l;\\
&{\rm if \ otherwise, \ } q_{i, l}' r_j = r_j q_{i, l}' \in r_j B_l \subset M_j B_l;\\
&{\rm in \  short, } \quad B_l r_j \subset M_j B_l;\\
&B_l B_j r_j= B_j B_l r_j \subset B_j M_j B_l \subset  M_j B_l;\\
&B_l r_l B_j r_j= B_j B_l r_l r_j= B_j B_l r_j r_l \subset B_j M_j
B_l r_l \subset M_j  B_l r_l.
\end {align*}
\end {proof}

\begin {lemma}
Consider the pair $(M, \tau)$ as described above and the GNS-
construction. Identify everything mentioned above as its image. We
$M''$ is the hyperfinite $II_1$ factor.
\end {lemma}

\begin{proof}
There is one and only one tracial state on $M_k$ for all $k \in
{\mathbb N}$. Hence the tracial state on $M$ is unique. Therefore
$M''$ is the hyperfinite $II_1$ factor, $R$.
\end {proof}

Define a unital *-endomorphism, $\Phi$, on $R$ to be the (right)
one-shift: i.e., sending $B_k$ to $B_{k+1}$, and sending $r_k$ to
$r_{k+1}$. We observe that $\Phi (R)$ is a $II_1$ factor and
\[
  [R: \Phi (R)]= n^2 < \infty.
\]

\begin {lemma}
The relative commutant $\Phi^k (R)' \cap R$ is exactly
$\otimes_{i=1}^k A$, on which the trace of $R$ is the product trace
given by the vector $\vec s$.
\end {lemma}

\begin {proof}
Because of our decomposition in Lemma 7 and Lemma 8, $R$ can be
written as
\[
  (\sum_{i=0}^j B_1 r_1^{i} B_1) \cdot (\sum_{i=0}^j
  B_2 r_2^{i} B_2) \cdot \cdots \cdot (\sum_{i=0}^j B_k r_k^{i} B_k)
  \cdot \Phi^k (R).
\]

Assume $x \in R \cap \Phi^k (R)'$. $x$ can be written, as in Lemma
8, of the following form:
\[
  x= \sum_{\vec {\alpha} \in \{ 0, \ 1,  \ \cdots, j-1 \}^k} y_1^{\vec {\alpha}} r_1^{g_1} z_1^{\vec
  {\alpha}}
  y_2^{\vec {\alpha}} r_2^{g_2} z_2^{\vec {\alpha}} \cdots y_k^{\vec {\alpha}} r_1^{g_k}
  z_k^{\vec {\alpha}}
  \cdot y^{\vec {\alpha}},
\]
where $\vec {\alpha}= (g_1, g_2, \dots, g_k)$ is a multi-index.
$y_1^{\vec {\alpha}}, z_1^{\vec {\alpha}} \in B_1$, $y_2^{\vec
{\alpha}}, z_2^{\vec {\alpha}} \in B_2$, $\cdots$,  $y_k^{\vec
{\alpha}}, z_k^{\vec {\alpha}} \in B_k$. $y^{\vec {\alpha}}$ is in
$\Phi^k (R)$. Note that $\Phi^k (R)$ is the weak closure of $\{
\Phi^k (M_i) \}_{i=1}^{\infty}$.

For every $\epsilon > 0$,  there exists an integer $i \in {\mathbb
N}$ such that
\begin{align*}
   &\forall \ {\vec {\alpha}}, \ z^{\vec {\alpha}} \in \Phi^k (M_i)
   \subset <B_{k+1}, r_{k+1}, \cdots, B_{k+i}, r_{k+i}>\\
   &\| x-\sum_{\vec {\alpha} \in \{ 0, \ 1,  \ \cdots, j-1 \}^k} y_1^{\vec {\alpha}} r_1^{g_1} z_1^{\vec
  {\alpha}}
  y_2^{\vec {\alpha}} r_2^{g_2} z_2^{\vec {\alpha}} \cdots y_k^{\vec {\alpha}} r_1^{g_k}
  z_k^{\vec {\alpha}}
  \cdot z^{\vec {\alpha}} \|_{2, \tau} < \delta\\
  &\delta= (\sqrt{\frac{j}{n}})^k \epsilon
\end{align*}

Put $L= l( l+1)/2+1$ for some integer $l > k+1$ and $l= 0 \mod 3$.
We have the following properties:
\begin {align*}
&[r_L, B_1]= [r_L, B_2]= \cdots = [r_L, B_{k+i}]= 0\\
&[r_L, r_2]= [r_L, r_3]= \cdots = [r_L, r_{k+i}]= 0\\
&r_L r_1 = \gamma r_1 r_L, \quad r_L r_1^{g_1}
   = \gamma^{g_1} r_1 r_L\\
&r_L r_L^*= r_L^* r_L= s_L\\
&r_L r_1 r_L^*= \gamma r_1 s_L, \quad r_L r_1^{g_1} r_L^*= \gamma^{g_1} r_1 s_L\\
&{\rm for \ } 0 \leq m \leq {j-1}, \quad r_L^m r_1^{g_1} (r_L^*)^m= \gamma^{g_1 m} r_1 s_l\\
&[s_L, B_1]= [s_L, B_2]= \cdots = [s_L, B_{k+i}]= 0\\
&[s_L, r_1]= [s_L, r_2]= \cdots = [s_L, r_{k+i}]= 0
\end {align*}

Therefore we claim:
\begin {align*}
&\| (x-\sum_{\vec {\alpha} \in \{ 0, \ 1,  \ \cdots, j-1 \}^k}
y_1^{\vec {\alpha}} r_1^{g_1} z_1^{\vec
  {\alpha}}
  y_2^{\vec {\alpha}} r_2^{g_2} z_2^{\vec {\alpha}} \cdots y_k^{\vec {\alpha}} r_1^{g_k}
  z_k^{\vec {\alpha}}
  \cdot z^{\vec {\alpha}}) s_L \|_{2, \tau}=\\
&\| (x-\sum_{\vec {\alpha} \in \{ 0, \ 1,  \ \cdots, j-1 \}^k}
y_1^{\vec {\alpha}} r_1^{g_1} z_1^{\vec
  {\alpha}}
  y_2^{\vec {\alpha}} r_2^{g_2} z_2^{\vec {\alpha}} \cdots y_k^{\vec {\alpha}} r_1^{g_k}
  z_k^{\vec {\alpha}}
  \cdot z^{\vec {\alpha}}) \frac{1}{j} \sum_{m=0}^{j-1} r_L^m {r_L^*}^m \|_{2, \tau}=\\
&\frac{1} {j} \| \sum_{\vec {\alpha}} \sum_m (r_L^m x {r_L^*}^m-
y_1^{\vec {\alpha}} r_L^m r_1^{g_1} (r_L^*)^m z_1^{\vec
  {\alpha}}
  y_2^{\vec {\alpha}} r_2^{g_2} z_2^{\vec {\alpha}} \cdots y_k^{\vec {\alpha}} r_1^{g_k}
  z_k^{\vec {\alpha}}
  \cdot z^{\vec {\alpha}}) \|_{2, \tau}=\\
&\frac{1} {j} \| \sum_{\vec {\alpha}} \sum_m (x - y_1^{\vec
{\alpha}}  r_1^{g_1 m} z_1^{\vec
  {\alpha}}
  y_2^{\vec {\alpha}} r_2^{g_2} z_2^{\vec {\alpha}} \cdots y_k^{\vec {\alpha}} r_1^{g_k}
  z_k^{\vec {\alpha}}
  \cdot z^{\vec {\alpha}}) s_L \|_{2, \tau}=\\
&\| (x-\sum_{\vec {\alpha} \in \{ 0, \ 1, \cdots, j-1 \}^k, g_1= 0}
y_1^{\vec {\alpha}} z_1^{\vec
  {\alpha}}
  y_2^{\vec {\alpha}} r_2^{g_2} z_2^{\vec {\alpha}} \cdots y_k^{\vec {\alpha}} r_1^{g_k}
  z_k^{\vec {\alpha}}
  \cdot z^{\vec {\alpha}}) s_L \|_{2, \tau}=\\
&\sqrt{\frac{j}{n}} \| x-\sum_{\vec {\alpha} \in \{ 0, \ 1, \cdots,
j-1 \}^k, g_1= 0} y_1^{\vec {\alpha}} z_1^{\vec
  {\alpha}}
  y_2^{\vec {\alpha}} r_2^{g_2} z_2^{\vec {\alpha}} \cdots y_k^{\vec {\alpha}} r_1^{g_k}
  z_k^{\vec {\alpha}}
  \cdot z^{\vec {\alpha}}  \|_{2, \tau}
\end{align*}

Since
\begin{align*}
&\{ x, \ y_1^{\vec {\alpha}}, \ z_1^{\vec
  {\alpha}}, \ y_2^{\vec {\alpha}}, \ r_2^{g_2}, \ z_2^{\vec {\alpha}}, \cdots,
 y_k^{\vec{\alpha}},
  \ r_1^{g_k}, \ z_k^{\vec {\alpha}}, \ z^{\vec {\alpha}} \} \subset \{ s_L, \ r_L, \ B_L
  \}'\\
&{\rm and \ } \tau (s_L)= \frac{j}{n}.
\end {align*}
Note that $\{ s_L, \ r_L, \ B_L \}'' = M_n({\mathbb C})_L$ is a type
$I$ factor \cite{sP83}.

By induction,
\begin{align*}
&\| x-\sum_{\vec {\alpha} \in \{ 0, \ 1, \cdots, j-1 \}^k, \ g_1= 0}
y_1^{\vec {\alpha}} z_1^{\vec
  {\alpha}}
  y_2^{\vec {\alpha}} r_2^{g_2} z_2^{\vec {\alpha}} \cdots y_k^{\vec {\alpha}} r_1^{g_k}
  z_k^{\vec {\alpha}}
  \cdot z^{\vec {\alpha}}  \|_{2, \tau} < \sqrt{\frac{n}{j}} \delta\\
&\| x-\sum_{\vec {\alpha} \in \{ 0, \ 1, \cdots, j-1 \}^k, \ g_1=
g_2= 0} y_1^{\vec {\alpha}} z_1^{\vec
  {\alpha}}
  y_2^{\vec {\alpha}} z_2^{\vec {\alpha}} \cdots y_k^{\vec {\alpha}} r_1^{g_k}
  z_k^{\vec {\alpha}}
  \cdot z^{\vec {\alpha}}  \|_{2, \tau} < (\sqrt{\frac{n}{j}})^2 \delta\\
&\cdots \\
&\| x-\sum_{\vec {\alpha} \in \{ 0 \}^k} y_1^{\vec {\alpha}}
z_1^{\vec
  {\alpha}}
  y_2^{\vec {\alpha}} z_2^{\vec {\alpha}} \cdots y_k^{\vec {\alpha}}
  z_k^{\vec {\alpha}}
  \cdot z^{\vec {\alpha}}  \|_{2, \tau} < (\sqrt{\frac{n}{j}})^k \delta= \epsilon\\
\end{align*}

Put $L_1= l_1( l_1+1)/2+1$ for some integer $l_1 > k+1$ and $l_1= 1
\mod 3$. We have the following properties:
\begin{align*}
   &U_{L_1}:= \sum_{m_1=0}^{j-1} q_{m_1, L_1}'\\
   &U_{L_1}^n= {\bf 1}\\
   &[U_{L_1}, \Phi (M_{k+i-1})]= 0
\end{align*}

Similarly, put $L_2= l_2( l_2+1)/2+1$ for some integer $l_2 > k+1$
and $l_2= 2 \mod 3$. We have the following properties:
\begin{align*}
   &U_{L_2}:= \sum_{m_2=0}^{j-1} q_{m_2, L_2}'\\
   &U_{L_2}^n= {\bf 1}\\
   &[U_{L_2}, \Phi (M_{k+i-1})]= 0
\end{align*}

Therefore
\begin{align*}
&\| x-\sum_{\vec {\alpha} \in \{ 0 \}^k} y_1^{\vec {\alpha}}
z_1^{\vec
  {\alpha}}
  y_2^{\vec {\alpha}} z_2^{\vec {\alpha}} \cdots y_k^{\vec {\alpha}}
  z_k^{\vec {\alpha}}
  \cdot z^{\vec {\alpha}}  \|_{2, \tau} =\\
&\| x-\sum_{\vec {\alpha} \in \{ 0 \}^k} y_1^{\vec {\alpha}}
z_1^{\vec
  {\alpha}}
  y_2^{\vec {\alpha}} z_2^{\vec {\alpha}} \cdots y_k^{\vec {\alpha}}
  z_k^{\vec {\alpha}}
  \cdot z^{\vec {\alpha}}  \frac {1}{n} \sum_{m_1=0}^{n-1} U_{L_1}^{m_1} {U_{L_1}^*}^{m_1} \|_{2,
  \tau}=\\
&\| x-\sum_{\vec {\alpha} \in \{ 0 \}^k} \frac {1}{n}
\sum_{m_1=0}^{n-1} U_{L_1}^{m_1} (y_1^{\vec {\alpha}} z_1^{\vec
  {\alpha}}) {U_{L_1}^*}^{m_1}
\cdots y_k^{\vec {\alpha}}
  z_k^{\vec {\alpha}}
  \cdot z^{\vec {\alpha}}  \|_{2, \tau}=\\
&\| x-\sum_{\vec {\alpha} \in \{ 0 \}^k} \frac {1}{n^2}
\sum_{m_1, m_2=0}^{n-1} U_{L_1}^{m_1} (y_1^{\vec {\alpha}} z_1^{\vec
  {\alpha}}) {U_{L_1}^*}^{m_1}
\cdots y_k^{\vec {\alpha}}
  z_k^{\vec {\alpha}}
  \cdot z^{\vec {\alpha}}  U_{L_2}^{m_2} {U_{L_2}^*}^{m_2}
  \|_{2, \tau}=\\
&\| x-\sum_{\vec {\alpha} \in \{ 0 \}^k} \frac {1}{n^2} \sum_{m_1,
m_2=0}^{n-1} U_{L_2}^{m_2} U_{L_1}^{m_1} (y_1^{\vec {\alpha}}
z_1^{\vec
  {\alpha}}) {U_{L_1}^*}^{m_1} {U_{L_2}^*}^{m_2}
\cdots y_k^{\vec {\alpha}}
  z_k^{\vec {\alpha}}
  \cdot z^{\vec {\alpha}}
  \|_{2, \tau}= \\
&\| x-\sum_{\vec {\alpha} \in \{ 0 \}^k} x_1^{\vec {\alpha}}
  y_2^{\vec {\alpha}} z_2^{\vec {\alpha}} \cdots y_k^{\vec {\alpha}}
  z_k^{\vec {\alpha}}
  \cdot z^{\vec {\alpha}}
  \|_{2, \tau}
\end{align*}

Observe that
\[
  x_1^{\vec {\alpha}}:= \frac {1}{n^2} \sum_{m_1, m_2=0}^{n-1} U_{L_2}^{m_2} U_{L_1}^{m_1} 
(y_1^{\vec {\alpha}} z_1^{\vec
  {\alpha}}) {U_{L_1}^*}^{m_1} {U_{L_2}^*}^{m_2}
\]
is the conditional expectation of $y_1^{\vec {\alpha}} z_1^{\vec
  {\alpha}} \in B_1$ onto $A_1$.

By induction,
\begin{align*}
&\| x-\sum_{\vec {\alpha} \in \{ 0 \}^k} x_1^{\vec {\alpha}}
  (y_2^{\vec {\alpha}} z_2^{\vec {\alpha}}) \cdots y_k^{\vec {\alpha}}
  z_k^{\vec {\alpha}}
  \cdot z^{\vec {\alpha}}
  \|_{2, \tau}=\\
&\| x-\sum_{\vec {\alpha} \in \{ 0 \}^k} x_1^{\vec {\alpha}}
  x_2^{\vec {\alpha}}  (y_3^{\vec {\alpha}}
  z_3^{\vec {\alpha}}) \cdots y_k^{\vec {\alpha}}
  z_k^{\vec {\alpha}}
  \cdot z^{\vec {\alpha}}
  \|_{2, \tau}=\\
&\cdots\\
&\| x-\sum_{\vec {\alpha} \in \{ 0 \}^k} x_1^{\vec {\alpha}}
  x_2^{\vec {\alpha}}  \cdots x_k^{\vec {\alpha}}
  \cdot z^{\vec {\alpha}}
  \|_{2, \tau} < \epsilon
\end{align*}
where $x_1^{\vec {\alpha}} \in A_1$, $x_2^{\vec {\alpha}} \in A_2$,
$\cdots$, $x_k^{\vec {\alpha}} \in A_k$.

Note that the von Neumann algebra  $\{ x, \ x_1^{\vec {\alpha}}, \
x_2^{\vec {\alpha}}, \ x_3^{\vec {\alpha}}, \cdots x_k^{\vec
{\alpha}} \}''$ commutes with $\Phi^k (M)$, which is a $II_1$
factor. Any element in the former von Neumann algebra has a scalar
conditional expectation onto $\Phi^k (M)$. In short, the former von
Neumann algebra and $\Phi^k (M)$ are mutually orthogonal.

According to \cite {sP83}, we have:
\begin{align*}
  &\| x-\sum_{\vec {\alpha} \in \{ 0 \}^k} x_1^{\vec {\alpha}}
  x_2^{\vec {\alpha}} x_3^{\vec {\alpha}} \cdots x_k^{\vec {\alpha}}
  \cdot \tau (z^{\vec {\alpha}})  \|_{2, \tau} < \epsilon\\
  &\sum_{\vec {\alpha} \in \{ 0 \}^k} x_1^{\vec {\alpha}}
  x_2^{\vec {\alpha}} x_3^{\vec {\alpha}} \cdots x_k^{\vec {\alpha}}
  \cdot \tau (z^{\vec {\alpha}})  \in A_1 \cdot A_2 \cdots A_k= \otimes^k
  A
\end{align*}

\end{proof}

\section {Discussion} \label {S:disc}

The construction of $\Phi$ depends on the choice of the
anticommutation set $S_1$. In the case of n-unitary shifts,
different choices of $S (\Psi; u)$ give uncountably many
nonconjugate shifts and at least a countably infinite family of
shifts that are pairwise not outer conjugate. Not to mention in
\cite {mC90} the existence of uncountably many non-outer-conjugate
nonbinary shifts, exploiting different 2-cocycles on the group $G=
\oplus_i^{\infty} \mathbb {Z}_2^{(i)}$. Each of the above has a
counterpart in our construction.

Given a finite dimensional $C^*$-algebra $A$ with ${\rm {rank}} (A)=
n$ and an n-unitary shift with the anticommutation set $S (\Psi; u)$
\cite {mC87} satisfies the following condition:
\[
  S( \Psi; u)= \{k_i \ | \
  k_{i+1} > k_i, \text {\ and \ } k_{i+2}- k_{i+1} > k_{i+1}- k_i,
  \ \forall
  i \in \mathbb {N} \}
\]

Denote by $M (\Psi; u)$ to be the ascending union of $M_k (\Psi;
u)=$
\[
  <A_1 (\Psi; u), r_1 (\Psi; u), A_2 (\Psi; u), r_2 (\Psi; u),
  \cdots, A_k (\Psi; u), r_k (\Psi; u)>
\]
with the trace $\tau (\Psi; u)$. We have
\[
   A \simeq A_1(\Psi; u) \simeq A_2(\Psi; u) \cdots
\]

Similarly, consider the pair $(M (\Psi; u), \tau (\Psi; u))$ as
described above and the GNS construction. Identify everything
mentioned above as its image. The weak closure $M (\Psi; u)''$ is
the hyperfinite $II_1$ factor, $R$.

Define a unital *-endomorphism, $\Psi$, on $R$ to be the (right)
one-shift: i.e., sending $A_k (\Psi; u)$ to $A_{k+1} (\Psi; u)$, and
sending $r_k (\Psi; u)$ to $r_{k+1} (\Psi; u)$. We observe that
$\Psi (R)$ is a $II_1$ factor and $[R: \Psi (R)]= n^2$.
\[
  \Psi^k (R)' \cap R= \otimes^k A.
\]

A well-known conjugacy invariant is the Connes-St{\o}rmer
entropy.
We estimate the entropy \cite {eS00} in the
following paragraph.

\begin {lemma}
$H (\Psi) \geq \ln n= \ln [R: \Psi (R)]$ no matter of the choice of
the anticommutation set $S (\Psi, u)$.
\end {lemma}

\begin {proof}
\begin {align*}
&H (\Psi)= \lim_{j \to \infty} H (M_j, \Psi) \geq \lim_{j \to
\infty}
H (\otimes^j A, \Psi)\\
&= \lim_{j \to \infty} \lim_{k \to \infty} \frac {1} {k}H (\otimes^j
A, \Psi (\otimes^j A),
\cdots, \Psi^{k-1} (\otimes^j A))\\
&= \lim_{j \to \infty} \lim_{k \to \infty} \frac {1} {k} H
(\otimes^{j +k -1} A) = \lim_{j \to \infty} \lim_{k \to \infty}
\frac {j+k-1} {k} H(A)\\
&= \lim_{j \to \infty} H(A)= H(A)
\end {align*}
\end {proof}

\begin {lemma}
$H (\Psi) \leq \ln n= \ln [R: \Psi (R)]$ no matter of the choice of
the anticommutation set $S (\Psi, u)$.
\end {lemma}

\begin {proof}
\begin {align*}
&H (\Psi)= \lim_{j \to \infty} H (M_j, \Psi) = \lim_{j \to \infty}
\lim_{k \to \infty} \frac {1} {k} H( M_j, \Psi (M_j), \cdots, \Psi^{k-1} (M_j))\\
& \leq \lim_{j \to \infty} \lim_{k \to \infty} \frac {1} {k} H(
M_{j+ k- 1})= \lim_{j \to \infty} \lim_{k \to \infty} \frac {1} {k}
H(\otimes^{j+k-1} M_n ({\mathbb C}))\\
&=\lim_{j \to \infty} \lim_{k \to \infty} \frac {1} {k} (j+k -1) \ln
n= \ln n
\end {align*}

\end {proof}

If $A$ is not trivial, $H(A)$ is non-zero.
Therefore, 
\[
  H( M_n ({\mathbb C})) \geq H (\Psi) \geq H(A) >0
\]
 no matter of the
choice of the anticommutation set $S (\Psi, u)$.

We conjecture that there is an endomorphism $\Psi$ that gives the entropy 
$H(A)$. Moreover, if $\bf A$ is a finite dimensional $C^*$-algebra with an
arbitrary trace vector (not necessarily of rational entries), we conjecture
that there is a similar result to Theorem 2, parallel with the 
Murray-von Neumann construction of the hyperfinite $II_1$ factor $R$.

\begin {acknowledgments}
I wish to express gratitude toward Professor Jing Yu for
his hospitality and support during my stay at the National 
Center of Theoretical Sciences in Taiwan. 
\end {acknowledgments}

\end {document}